\input amstex
\documentstyle{amsppt}
\NoBlackBoxes
\input cyracc.def
\newfam\cyrfam
\catcode `\@=11
\font\tencyr=wncyr10
\font\eightcyr=wncyr8
\catcode `\@=13
\addto\tenpoint{\def\cyr{\fam\cyrfam\tencyr\cyracc}}
\addto\eightpoint{\def\cyr{\fam\cyrfam\eightcyr\cyracc}}
\textfont\cyrfam=\tencyr
\define\sh{{\cyr sh}}
\define\al{\alpha}
\define\be{\beta}
\define\de{\delta}
\define\ep{\epsilon}
\define\om{\omega}
\define\si{\sigma}
\define\zt{\zeta}
\define\De{\Delta}
\define\Si{\Sigma}

\define\HH{\frak H}

\define\QS{\operatorname{QSym}}

\define\End{\operatorname{End}}
\define\id{\operatorname{id}}
\define\QA{\bold Q\langle x,y\rangle}

\define\BJOP{1}
\define\BBB{2}
\define\BBBL{3}
\define\BK{4}
\define\Eh{5}
\define\Gei{6}
\define\Ges{7}
\define\G{8}
\define\HP{9}
\define\Ha{10}
\define\Hb{11}
\define\Hc{12}
\define\IK{13}
\define\K{14}
\define\LM{15}
\define\M{16}
\define\Oh{17}
\define\OZ{18}
\define\Ree{19}
\define\RSS{20}
\define\Z{21}
\topmatter
\title Relations of Multiple Zeta Values and Their Algebraic Expression
\endtitle
\author Michael E. Hoffman$^1$ and Yasuo Ohno$^2$\endauthor
\affil ${}^1$U. S. Naval Academy\\
${}^2$Kinki University\endaffil
\address Mathematics Department, U. S. Naval Academy, Annapolis, MD 21402 
USA\endaddress
\email meh\@usna.edu\endemail
\address Department of Mathematics and Physics, Faculty of Science and 
Technology, Kinki University, 3-4-1, Kowakae, Higashi-Osaka, Osaka, 
577-8502 Japan\endaddress
\email ohno\@math.kindai.ac.jp\endemail
\abstract
We establish a new class of relations among the multiple zeta values
$$
\zt(k_1,\dots,k_l)=\sum_{n_1>\dots>n_l\ge 1}\frac1{n_1^{k_1}\cdots n_k^{k_l}} ,
$$
which we call the cyclic sum identities.  These identities have an
elementary proof, and imply the ``sum theorem'' for multiple zeta
values.  They also have a succinct statement
in terms of ``cyclic derivations'' as introduced by Rota, Sagan and
Stein.  In addition, we discuss the expression of other relations of multiple 
zeta values via the shuffle and ``harmonic'' products on the underlying vector 
space $\HH$ the noncommutative polynomial ring $\QA$, and also 
using an action of the Hopf algebra of quasi-symmetric functions on $\QA$.
\endabstract
\leftheadtext{Michael E. Hoffman and Yasuo Ohno}
\rightheadtext{Relations of multiple zeta values}
\endtopmatter
\document
\par
\subheading{1. Multiple zeta values}
Let $k_1, k_2,\dots, k_l$ be positive integers with $k_1>1$.  The
multiple zeta value $\zt(k_1,k_2\dots,k_l)$ (of weight $k_1+\dots+k_l$
and length $l$) associated with this sequence is 
the sum of the convergent $l$-fold infinite series
$$
\sum_{n_1>n_2>\cdots>n_l\ge 1}\frac1{n_1^{k_1}n_2^{k_2}\cdots n_l^{k_l}} .
$$
These quantities were introduced under the name ``multiple harmonic series''
in \cite{\Ha}, and independently (with an opposite convention for
the order of the sequence) in \cite{\Z}.  
They have appeared in knot theory \cite{\LM}, quantum field theory
\cite{\BK}, and in connection with mirror symmetry \cite{\Hc}.
There are many relations among multiple zeta values (henceforth MZVs), 
starting with
the widely known and often rediscovered identity $\zt(2,1)=\zt(3)$.
An outstanding example is the ``sum theorem'', which
says that the sum of all MZVs of fixed length and weight is independent 
of length, i.e.
$$
\sum_{\{(k_1,\dots,k_l)|k_1+\dots+k_l=n,k_1>1\}}\zt(k_1,k_2,\dots,k_l)=\zt(n) .
$$
This was conjectured in \cite{\Ha} and proved independently by 
A. Granville \cite{\G} and D. Zagier.  
Many other identities have been conjectured and proved in the last
decade \cite{\BBB,\BBBL,\LM,\Oh,\OZ}, but surprising new ones continue 
to appear.  
In this paper we establish a new class of relations of MZVs,
which can be stated as follows.
\proclaim{Cyclic sum theorem} For any positive integers $k_1,k_2,\dots,k_l$
with some $k_i\ge 2$,
$$\multline
\sum_{j=1}^l \zt(k_j+1,k_{j+1},\dots,k_l,k_1,\dots,k_{j-1})=\\
\sum_{\{j|k_j\ge2\}}\sum_{q=0}^{k_j-2}\zt(k_j-q,k_{j+1},\dots,k_l,k_1,\dots,
k_{j-1},q+1) .
\endmultline\tag1$$
\endproclaim
\par
The name we have given this result will be clearer if we state it in
an alternative form using the ``duality'' of MZVs.
If $s=(k_1,\dots,k_l)$ is an admissible sequence of positive integers 
(i.e., $k_1>1$) with sum $n$, its dual sequence $\tau(s)=(j_1,\dots,j_{n-l})$ 
(having sum $n$ and $j_1>1$) is defined as follows.  Let $\Si$ be the map 
that takes a sequence to its sequence of partial sums, and let $\frak I_n$ 
be the set of strictly increasing
sequences of positive integers with last element at most $n$.
Set $\tau(s)=\Si^{-1}C_nR_n\Si(s)$, where for $(a_1,\dots,a_i)\in\frak I_n$,
$$
R_n(a_1,\dots,a_i)=(n+1-a_i,n+1-a_{i-1},\dots,n+1-a_1)
$$
and $C_n(a_1,\dots,a_i)$ is the complement of $\{a_1,\dots,a_i\}$ in
$\{1,2,\dots,n\}$ arranged in increasing order.
Then $\tau$ is an involution on the set of admissible 
positive-integer sequences.
Now call two sequences of positive integers 
cyclically equivalent if one is a cyclic permutation of the other, and
let $\Pi(n,l)$ be the set of cyclic equivalence classes of (not necessarily
admissible) positive-integer sequences of sum $n$ and length $l$.  
For the juxtaposition $s_1s_2$ of two admissible sequences one has 
$\tau(s_1s_2)=\tau(s_2)\tau(s_1)$ (see \cite\Ha, Proposition 3.1), so
$\tau(s')$ is a cyclic permutation of $\tau(s)$ when $s'$ is an (admissible) 
cyclic permutation of $s$, and thus any equivalence class $[s]\in\Pi(n,l)$ 
has a dual equivalence class $[\tau(s)]\in\Pi(n,n-l)$.
(For example, $\{(2,3),(3,2)\}\in\Pi(5,2)$ has dual equivalence class
$\{(2,1,2),(1,2,2),(2,2,1)\}\in\Pi(5,3)$.)
The cyclic sum theorem can be restated as follows:  for any admissible sequence
$s=(k_1,\dots,k_l)$ with $k_1+\dots+k_l=n$,
$$
\sum_{(p_1,\dots,p_l)\in [s]}\zt(p_1+1,p_2,\dots,p_l)=
\sum_{(q_1,\dots,q_{n-l})\in [\tau(s)]}\zt(q_1+1,q_2,\dots,q_{n-l}).
\tag2$$
(For the equivalence of the two forms, see the remarks following 
Theorem 2.3 below.)
\par
The cyclic sum theorem has an elementary proof involving partial fractions,
but admits a remarkably simple expression in terms of ``cyclic derivations''
(in the sense of \cite{\RSS}) of the noncommutative polynomial algebra $\QA$.  
In this way it parallels an earlier result (Theorem 5.1 of
\cite{\Ha}, reformulated as Theorem 2.1 below), which 
was proved by an unenlightening partial-fractions argument but
can be expressed very simply in terms of ordinary derivations of $\QA$.  
The cyclic sum theorem also implies the sum theorem, giving a
new proof of this result which does not involve generating functions
as used by Granville and Zagier.
\par
We introduce our algebraic machinery in \S2: as in \cite{\Hb}, we 
think of MZVs as images of monomials under a map
$\zt:\HH^0\to\bold R$, where $\HH^0$ is an appropriate subspace of
$\QA$.  We define derivations and cyclic derivations of $\QA$
that give rise to identities of MZVs, and show how
the cyclic sum theorem implies the sum theorem.  In \S3 we prove the
cyclic sum theorem by elementary methods.  In \S4 we return to 
our algebraic viewpoint, 
recalling the shuffle and ``harmonic'' products on $\HH^0$
(both of which make $\zt$ a homomorphism), and relating them to
Theorem 2.1.  In \S5 we introduce an 
action of the Hopf algebra $\QS$ of quasi-symmetric functions
on $\QA$.  We show how Theorem 2.1, the cyclic sum theorem, and 
the identities proved by the second author in \cite{\Oh} can be 
expressed in terms of this action.
\par
The first author conjectured the cyclic sum theorem in August 1999,
and thanks Michael Bigotte for checking it by computer against
tables of known relations \cite{\BJOP} through weight 12.
The second
author proved the conjecture during his stay at the Max-Planck-Institut
f\"ur Mathematik in Bonn in early 2000, and he thanks Masanobu Kaneko
and Don Zagier for useful discussions and the Institut for its hospitality.
The second author is supported in part by the Research Fellowship of 
the Japan Society for the Promotion of Science for Young Scientists.
\subheading{2. A noncommutative polynomial algebra and its derivations}  
In this section we introduce an algebraic approach by thinking of MZVs
as values of a homomorphism from a subspace of the noncommutative
polynomial algebra $\QA$ to the reals.  We then consider both derivations
and ``cyclic derivations'' (defined below) of $\QA$, and formulate 
relations of MZVs in terms of them.  We state the cyclic sum theorem
algebraically (Theorem 2.3), show it is equivalent to the two forms 
(1) and (2) given in \S1, and finally prove the sum theorem from the 
cyclic sum theorem.
\par
Let $\QA$ be the algebra of polynomials over the rationals
in noncommutative indeterminates $x,y$, regarded as a graded 
$\bold Q$-algebra with $x$ and $y$ both of degree 1.  For any word
(monomial) $w$ of $\QA$, denote by $|w|$ its total degree (also called
its weight) and by $\ell(w)$ the number of occurrences of $y$ in $w$
(called the length of $w$).  We call $|w|-\ell(w)$ the colength of
$w$:  it is the number of occurrences of $x$ in $w$.  
The underlying graded rational vector space of $\QA$ is denoted $\HH$.
\par
Let $\HH^1=\bold Q1\oplus\HH y$ and $\HH^0=\bold Q 1\oplus x\HH y$.
Then $\HH^1$ is a subalgebra of $\QA$, in fact the noncommutative 
polynomial algebra on generators $z_i=x^{i-1}y$.
We have $\HH=\HH x\oplus\HH^1$ and $\HH^1=y\HH^1\oplus\HH^0$.
We can think of MZVs as images of words of $\QA$
under the $\bold Q$-linear map $\zt: \HH^0\to\bold R$ defined by
$\zt(1)=1$ and
$$
\zt(x^{k_1-1}yx^{k_2-1}y\cdots x^{k_l-1}y)=\zt(k_1,k_2,\dots,k_l)
$$
for any positive integers $k_1,k_2,\dots,k_l$ with $k_1>1$.
\par
Let $\tau$ be the anti-automorphism of $\QA$ exchanging $x$ 
and $y$, e.g. $\tau(x^2yxy)=xyxy^2$.  Evidently $\tau$ is an involution.
Applied to words, $\tau$ preserves weight and exchanges
length and colength: note that $\HH^0$ (but not $\HH^1$) is closed 
under $\tau$.  It is easy to check that for dual sequences
$s=(k_1,\dots,k_l)$ and $\tau(s)=(j_1,\dots,j_{n-l})$ as defined in
the preceding section, $\tau(x^{k_1-1}y\cdots x^{k_l-1}y)=
x^{j_1-1}y\cdots x^{j_{n-l}-1}y$.
\par
As usual, by a derivation of $\QA$ we mean a map $F:\HH\to\HH$ (of
graded rational vector spaces) such that $F(uv)=F(u)v+uF(v)$ for
all $u,v\in\HH$.  The
commutator of two derivations is a derivation, so the set
of derivations of $\QA$ is a Lie algebra graded by degree.  
If $\de$ is a derivation, then $\bar\de=\tau\de\tau$
is also a derivation (of the same degree).  We call a derivation $\de$
symmetric if $\bar\de=\de$ and antisymmetric if $\bar\de=-\de$.
Since $\overline{[\de,\ep]}=[\bar\de,\bar\ep]$, symmetric and antisymmetric
derivations behave nicely under commutator, e.g. 
$[\de,\ep]$ is antisymmetric
if $\de$ is symmetric and $\ep$ antisymmetric.  Note that a symmetric
or antisymmetric derivation is completely determined by where it sends
$x$.
\par
We denote by $D$ the derivation such that $D(x)=0$ and $D(y)=xy$.
In terms of the generators $z_i$ of $\HH^1$
mentioned above, we have $D(z_i)=z_{i+1}$ and more generally
$$
D(z_{i_1}z_{i_2}\cdots z_{i_l})=
z_{i_1+1}z_{i_2}\cdots z_{i_l}+
z_{i_1}z_{i_2+1}z_{i_3}\cdots z_{i_l}+\dots+
z_{i_1}\cdots z_{i_{l-1}}z_{i_l+1} .
\tag3$$
The following result was proved by a partial-fractions argument in
\cite{\Ha}.  We note that the hypothesis on $w$ cannot be weakened,
since $\zt(D(y))=\zt(xy)\ne 0=\zt(\bar D(y))$.
\proclaim{Theorem 2.1} For any word $w$ of $\HH^0$, $\zt(D(w))
=\zt(\bar D(w))$.
\endproclaim
K. Ihara and M. Kaneko \cite{\IK} have generalized Theorem 2.1
as follows.  For $n\ge 1$ let $\partial_n$
be the antisymmetric derivation with $\partial_n(x)=x(x+y)^{n-1}y$;
note that $\partial_1=\bar D-D$.
\proclaim{Theorem (Ihara and Kaneko)} For all $n\ge 1$ and words $w$ of
$\HH^0$, $\zt(\partial_n(w))=0$.
\endproclaim
We shall discuss the proof of this result in \S5 below.
\par
We now consider cyclic derivations of $\QA$; these are not derivations,
but rather are defined as follows (cf. \cite{\RSS}).
\proclaim{Definition}
A cyclic derivation $\psi$ of $\QA$ is a $\bold Q$-linear map $\psi:\HH\to
\End\HH$, where $\End\HH$ is the graded rational vector space of 
endomorphisms of $\HH$ (as a graded rational vector space), such that
$$
(\psi(f_1f_2),f)=(\psi(f_1),f_2f)+(\psi(f_2),ff_1)
\tag4$$
for all $f_1,f_2,f\in\HH$, where $(\be,u)$ denotes the image of $u\in\HH$
under $\be\in\End\HH$.
\endproclaim
If $\psi$ is a cyclic derivation, evidently $\psi(1)$ is the zero endomorphism.
By induction equation (4) is easily extended to the identity
$$\multline
(\psi(f_1f_2\cdots f_n),f)=(\psi(f_1),f_2\cdots f_nf)+
(\psi(f_2),f_3\cdots f_nff_1)\\
+\dots+(\psi(f_n),ff_1\cdots f_{n-1})
\endmultline\tag5$$
(cf. Proposition 2.3 of \cite{\RSS}; unfortunately the statement has a
misprint).
\par
Just as the conjugate by $\tau$ of an ordinary derivation is a derivation,
it is possible to conjugate a cyclic derivation by $\tau$ as follows.
\proclaim{Proposition 2.2} Suppose $\psi$ is a cyclic derivation.  Then
the map $\bar\psi:\HH\to\End\HH$ is also a cyclic derivation, where
$$
(\bar\psi(f),g)=\tau(\psi(\tau(f)),\tau(g)) 
$$
for $f,g\in\HH$.
\endproclaim
\demo{Proof} It suffices to check identity (4), which is routine.
\qed\enddemo
\par
If $\psi$ is a cyclic derivation and $f\in\HH$, the endomorphism
$\psi(f)$ gives rise to a canonical element $(\psi(f),1)$ of $\HH$;
we shall abuse notation and write $\psi(f)$ for this element of $\HH$
when no confusion can arise.  Note that as elements of $\HH$, 
$\psi(fg)=\psi(gf)$ and $\bar\psi(f)=\tau\psi(\tau(f))$ for any
$f,g\in\HH$.
\par
Now we define the cyclic derivation $C$ of $\QA$ by setting $C(x)=0$
(zero endomorphism) and $(C(y),f)=xfy$ for all $f\in\HH$.  By applying
identity (5), it is easy to see that
$$
(C(z_i),f)=(C(x^{i-1}y),f)=(C(y),fx^{i-1})=xfx^{i-1}y=xfz_i
$$
for $f\in\HH$, and thus (abusing notation as indicated above) that
$C(z_i)=z_{i+1}$.  Given an arbitrary monomial $w=z_{i_1}z_{i_2}
\cdots z_{i_l}$ of $\HH^1$, we can apply identity (5) to get
$$
C(w)=z_{i_1+1}z_{i_2}\cdots z_{i_l}+
z_{i_2+1}z_{i_3}\cdots z_{i_l}z_{i_1}+\dots+
z_{i_l+1}z_{i_1}\cdots z_{i_{l-1}} ,
\tag6
$$
which may be compared to equation (3).
Similarly, it can be shown that
$$
\bar C(z_{i_1}\cdots z_{i_l})=\sum_{i_j\ge 2}\sum_{q=0}^{i_j-2}z_{i_j-q}
z_{i_{j+1}}\cdots z_{i_l}z_{i_1}\cdots z_{i_{j-1}}z_{q+1} .
$$
Thus, form (1) of cyclic sum theorem as stated in \S1 may be expressed as
follows; the proof is given in \S3 below.
\proclaim{Theorem 2.3} For any word $w$ of $\HH^1$ that is not a power of 
$y$, $\zt(C(w))=\zt(\bar C(w))$.
\endproclaim
Let $s=(k_1,\dots,k_l)$ be an admissible sequence of positive integers,
and $w=x^{k_1-1}y\cdots x^{k_l-1}y$ the corresponding word of $\HH^0$.  
In view of equation (6),
$$
\zt(C(w))=m(w)\sum_{(p_1,\dots,p_l)\in [s]}\zt(p_1+1,p_2,\dots,p_l),
$$
where $[s]$ is the equivalence class of $s$ in $\Pi(k_1+\dots+k_l,l)$
and $m(w)$ is the largest integer $m$ such that $w=u^m$ for $u\in\HH^0$.
Then form (2) of the cyclic sum theorem is equivalent to 
$\zt(C(w))=\zt(C\tau(w))$ (note $m(\tau(w))=m(w)$), and
this is equivalent to Theorem 2.3 since $\zt$ is $\tau$-invariant 
(see Theorem 4.1 below).
\par
Theorems 2.1 and 2.3 are of course formally very similar.  Both give
an equation between a sum of MZVs of length $l$ and a sum of MZVs
of length $l+1$.  
An important
difference between the two results is that $C$ is much simpler than $D$
on periodic words of $\HH^0$.  For example, Theorem 2.3 applied to $z_n^l$ 
gives
$$
\zt(z_{n+1}z_n^{l-1})=\sum_{i=0}^{n-2}\zt(z_{n-i}z_n^{l-1}z_{i+1}) .
$$
In the sequence notation, the case $n=3$ is
$$
\zt(4,3,\dots,3)=\zt(3,3,\dots,3,1)+\zt(2,3,\dots,3,2) ,
$$
which does not seem to follow easily from other known identities.
We close this section by deducing the sum theorem from the cyclic sum
theorem.
\proclaim{Corollary 2.4} For any integers $1\le l<n$, let $S(n,l)$ be
the sum of words $w\in\HH^0$ with $|w|=n$ and $\ell(w)=l$.  Then 
$\sum_{w\in S(n,l)} \zt(w)$ is independent of $l$ (and, in particular,
is equal to $\zt(z_n)=\zt(n)$).
\endproclaim
\demo{Proof} Consider the element $\mu=(x+ty)^{n-1}-x^{n-1}-t^{n-1}y^{n-1}
\in\HH^0[t]$.  From the properties of cyclic derivation (cf. Corollary 2.5
of \cite{\RSS}) we have
$$\multline
C((x+ty)^{n-1})=(C((x+ty)^{n-1}),1)=(n-1)(C(x+ty),(x+ty)^{n-2})\\
=(n-1)(tC(y),(x+ty)^{n-2})=(n-1)tx(x+ty)^{n-2}y,
\endmultline$$
and thus $C(\mu)=(n-1)(tx(x+ty)^{n-2}y-t^{n-1}xy^{n-1})$; a similar
calculation gives $\bar C(\mu)=(n-1)(x(x+ty)^{n-2}y-x^{n-1}y)$.  For each
$1\le l<n-1$, the coefficient of $t^l$ in $\zt(\bar C(\mu)-C(\mu))$ gives
the identity
$$
\sum_{w\in S(n,l)}\zt(w)-\sum_{w\in S(n,l+1)}\zt(w)=0
$$
after dividing by $n-1$.
\qed\enddemo
\par
\subheading{3. Proof of the cyclic sum theorem}  For positive
integers $k_1,k_2,\dots, k_l$ and nonnegative integer $k_{l+1}$,
let
$$
T(k_1,\dots,k_l)=\sum_{n_1>n_2>\dots>n_l>n_{l+1}\ge 0}
\frac1{(n_1-n_{l+1})n_1^{k_1}\dots n_l^{k_l}}
$$
and 
$$
S(k_1,\dots,k_l,k_{l+1})=\sum_{n_1>n_2>\dots>n_l>n_{l+1}> 0}
\frac1{(n_1-n_{l+1})n_1^{k_1}\dots n_l^{k_l}n_{l+1}^{k_{l+1}}} .
$$
For the convergence of these series, we have the following.
\proclaim{Theorem 3.1} $T(k_1,\dots,k_l)$ is bounded when one of 
$k_1,\dots,k_l$ exceeds 1, and $S(k_1,\dots,k_l,k_{l+1})$ is 
bounded when one of $k_1,\dots,k_l,k_{l+1}+1$ exceeds 1.
\endproclaim
Our key result is as follows.
\proclaim{Theorem 3.2} For any positive integers $k_1,k_2\dots,k_l$
with $k_i>1$ for some $i$,
$$
T(k_1,\dots,k_l)-T(k_2,\dots,k_l,k_1)=\zt(k_1+1,k_2,\dots,k_l)-\sum_{j=0}^{k_1-2}
\zt(k_1-j,k_2,\dots,k_l,j+1)
$$
where the sum on the right is understood as 0 if $k_1=1$.
\endproclaim
\par
To prove the cyclic sum theorem in the form stated in \S1, sum
Theorem 3.2 over all cyclic permutations of the sequence $(k_1,\dots,k_l)$.
\par
It is immediate that
$$
S(k_1,\dots,k_l,0)=T(k_1,\dots,k_l)-\zt(k_1+1,k_2,\dots,k_l) .
\tag7$$
Also, applying the identity
$$
\frac1{n_1(n_1-n_{l+1})}=\frac1{n_{l+1}}\left(\frac1{n_1-n_{l+1}}-\frac1{n_1}
\right)
\tag8$$
to $S(k_1,\dots,k_l,k_{l+1})$ gives
$$
S(k_1,\dots,k_l,k_{l+1})=S(k_1-1,k_2,\dots,k_l,k_{l+1}+1)-\zt(k_1,\dots,k_l,
k_{l+1}+1) .
\tag9$$
Finally, applying (8) to $S(1,k_2,\dots,k_l,k_{l+1})$ gives
$$\align
&\sum_{n_1>n_2>\dots>n_l>n_{l+1}> 0}
\frac1{n_2^{k_2}\dots n_l^{k_l}n_{l+1}^{k_{l+1}+1}} 
\left(\frac1{n_1-n_{l+1}}-\frac1{n_1}\right)= \\
&\sum_{n_2>\dots>n_l>n_{l+1}> 0}
\frac1{n_2^{k_2}\dots n_l^{k_l}n_{l+1}^{k_{l+1}+1}}
\sum_{n_1=n_2+1}^{\infty}\left(\frac1{n_1-n_{l+1}}-\frac1{n_1}\right)= \\
&\sum_{n_2>\dots>n_l>n_{l+1}> 0}
\frac1{n_2^{k_2}\dots n_l^{k_l}n_{l+1}^{k_{l+1}+1}}
\sum_{j=0}^{n_{l+1}-1}\frac1{n_2-j}= \\
&\sum_{n_2>\dots>n_l>n_{l+1}>j\ge 0}
\frac1{(n_2-j)n_2^{k_2}\dots n_l^{k_l}n_{l+1}^{k_{l+1}+1}}
\endalign$$
and so
$$
S(1,k_2,\dots,k_l,k_{l+1})=T(k_2,\dots,k_l,k_{l+1}+1) .
\tag{10}$$
\par
\demo{Proof of Theorem 3.2} Apply equation (7), then equation (9) 
$k_1-1$ times, and finally equation (10):
$$\multline
T(k_1,\dots,k_l)-\zt(k_1+1,k_2,\dots,k_l)=S(k_1,\dots,k_l,0)=\\
S(k_1-1,k_2,\dots,k_l,1)-\zt(k_1,\dots,k_l,1)=\dots=\\
S(1,k_2,\dots,k_l,k_1-1)
-\sum_{j=0}^{k_1-2}\zt(k_1-j,k_2,\dots,k_l,j+1)=\\
T(k_2,\dots,k_l,k_1)
-\sum_{j=0}^{k_1-2}\zt(k_1-j,k_2,\dots,k_l,j+1).
\qed\endmultline$$
\enddemo
\demo{Proof of Theorem 3.1} Using equation (7),
$$
S(k_1,\dots,k_l,k_{l+1})\le S(k_1,\dots,k_l,0)\le T(k_1,\dots,k_l),
$$
so $S(k_1,\dots,k_{l+1})$ is bounded if
$T(k_1,\dots,k_l)$ is; and if $k_1=1$, equation (10) says
$S(k_1,\dots,k_{l+1})=T(k_2,\dots,k_l,k_{l+1}+1)$.
So the statement about the $S$'s 
follows from the one about the $T$'s.  Also, to prove the first
assertion it is evidently enough to treat the case $k_1+\dots+k_l=l+1$.  
Now 
$$\align
T(2,1,\dots,1)&=\sum_{n_1>n_2>\dots>n_{l+1}\ge 0}
\frac1{n_1^2(n_1-n_{l+1})n_2\dots n_l} \\
&\le\sum\Sb n_1>n_2>\dots>n_l>0\\n_1\ge j>0\endSb\frac1{n_1^2jn_2\dots n_l}\\
&=\zt(3,\underbrace{1,\dots,1}_{l-1})+l\zt(2,\underbrace{1,\dots,1}_l)+
\sum_{i=1}^{l-1}\zt(2,\underbrace{1,\dots,1}_{i-1},2,
\underbrace{1,\dots,1}_{l-i-1}),
\endalign$$
so $T(2,1,\dots,1)$ is bounded.  Then by equations (7) and (10), we have
$$\align
T(1,2,1,\dots,1)&=S(1,2,1,\dots,1,0)+\zt(2,2,1,\dots,1)\\
&=T(2,1,\dots,1,1)+\zt(2,2,1,\dots,1)
\endalign$$
and we can continue in this way to bound all the sums
$T(1,\dots,1,2,1,\dots,1)$. \qed
\enddemo
\subheading{4. Commutative multiplications on $\HH$} There are two
commutative multiplications on the vector space $\HH$, both of which
have significance for MZVs.  First, there is the
shuffle product $\sh$, which can be defined inductively on words of $\HH$
by requiring that it distribute over addition and satisfy the axioms
\roster
\item"{S1.}"
for any word $w$, $1\sh w=w\sh 1=w$;
\item"{S2.}"
for any words $w_1,w_2$ and $a,b\in\{x,y\}$, 
$$
aw_1\sh bw_2=a(w_1\sh bw_2)+b(aw_1\sh w_2) .
$$
\endroster
It is evident that $(\HH^0,\sh)$ is a subalgebra of $(\HH,\sh)$, and
we have the following result.
\proclaim{Theorem 4.1} $\zt$ is a $\tau$-invariant homomorphism of 
$(\HH^0,\sh)$ into $\bold R$.
\endproclaim
\demo{Proof}  This follows from the representation of MZVs
as iterated integrals (see \cite{\Z,\K,\Hb,\HP}).  If we define iterated integrals
recursively by
$$
\int_0^t\al_1=\int_0^t f(s)ds
$$
and 
$$
\int_0^t \al_1\al_2\cdots \al_n=\int_0^t f(s)\left(\int_0^s \al_2\cdots\al_n
\right)ds
$$
for $\al_1=f(t)dt$, then it is easy to show that for any nonnegative
integers $p_1,p_2,\dots,p_l$ with $p_1>1$,
$$
\zt(x^{p_1}y\cdots x^{p_l}y)=\int_0^1\om_0^{p_1}\om_1\cdots\om_0^{p_l}\om_1 
$$
where $\om_0=dt/t$ and $\om_1=dt/(1-t)$.  That $\zt$ is a homomorphism
then follows from the fact that iterated integrals multiply according
to shuffle products \cite{\Ree}.  The $\tau$-invariance follows from
a change of variable.
\qed\enddemo
Second, there is the ``harmonic'' product $*$ defined on $\HH$ 
by requiring that it distribute over addition and satisfy the axioms
\roster
\item"{H1.}"
for any word $w$, $1*w=w*1=w$;
\item"{H2.}"
for any word $w$ and positive integer $p$, $x^p*w=wx^p$;
\item"{H3.}"
for any words $w_1, w_2$ and positive integers $p,q$,
$$
x^{p-1}yw_1*x^{q-1}yw_2=x^{p-1}y(w_1*x^{q-1}yw_2)+x^{q-1}y(x^{p-1}yw_1*w_2)
+x^{p+q-1}y(w_1*w_2).
$$
\endroster
As was shown in \cite{\Hb}, this defines a commutative and associative
product on $\HH$.  If $\HH^1$ is regarded as the underlying vector space
of the noncommutative algebra $\bold Q\langle z_1,z_2,\dots\rangle$,
where $z_i=x^{i-1}y$ as in \S2, then axiom (H3) for words of $\HH^1$ reads
$$
z_pw_1*z_qw_2=z_p(w_1*z_qw_2)+z_q(z_pw_1*w_2)+z_{p+q}(w_1*w_2),
$$
which may be compared to (S2).  This can be thought of as describing
multiplication of series, and the following result is proved in \cite{\Hb}.
\par
\proclaim{Theorem 4.2} The map $\zt:(\HH^0,*)\to\bold R$ is a homomorphism.
\endproclaim
\par
On the other hand, we can define a $\bold Q$-linear map from 
$\phi:\HH^1\to\bold Q[[t_1,t_2,\dots]]$, where 
$\bold Q[[t_1,t_2,\dots]]$ is the $\bold Q$-algebra of 
formal series in the countable set of (commuting) variables $t_1,t_2,\dots$,
by setting $\phi(1)=1$ and 
$$
\phi(z_{i_1}\cdots z_{i_l})=\sum_{n_1>n_2>\dots>n_l\ge 1}
t_{n_1}^{i_1}t_{n_2}^{i_2}\cdots t_{n_l}^{i_l} .
$$
Then $\phi$ is a homomorphism, and in fact a monomorphism; its image
is the algebra $\QS$ of quasi-symmetric functions 
as defined in \cite{\Ges}.  A formal power series (of bounded degree)
in $t_1,t_2,\dots$ is called a quasi-symmetric function
if the coefficients and $t_{i_1}^{p_1}t_{i_2}^{p_2}\cdots t_{i_k}^{p_k}$
and $t_{j_1}^{p_1}t_{j_2}^{p_2}\cdots t_{j_k}^{p_k}$ are the same whenever
$i_1<i_2<\cdots<i_k$ and $j_1<j_2<\cdots<j_k$.  Evidently any 
quasi-symmetric function is a sum of monomial quasi-symmetric functions,
which can be defined as the elements $\phi(z_{i_1}\cdots z_{i_k})$.
\par
Since $\zt:\HH^0\to\bold R$
is a homomorphism for both multiplications, any element of the form
$\langle u,v\rangle=u\sh v - u*v$ for $u,v\in\HH^0$ must be in the 
kernel of $\zt$.
Together with Theorem 2.1, our next result shows that this remains
true for elements of the form $\langle y,w\rangle$, $w\in\HH^0$.
\proclaim{Theorem 4.3} For words $w$ of $\HH$, $y\sh w-y*w=\bar D(w)-D(w)$,
where $D$ is the derivation of \S2.
\endproclaim
\demo{Proof} We proceed by induction on $\ell(w)$.  If $\ell(w)=0$, then
$w=x^k$ and we have
$$
y\sh x^k-y*x^k=yx^k+xyx^{k-1}+\dots+x^ky-yx^k=
xyx^{k-1}+\dots+x^ky=\bar D(x^k)- D(x^k).
$$
Now suppose $\ell(w)=n>0$ and the result holds for words of length
less than $n$.  Then we can write $w=x^kyw_1$ for some word $w_1$
with $\ell(w_1)=n-1$.  Using the inductive definitions (S2) and (H3),
we have
$$\align
y\sh(x^kyw_1)&=yx^kyw_1+xyx^{k-1}yw_1+\dots+x^ky^2w_1+x^ky(y\sh w_1)\\
y*(x^kyw_1)&=yx^kyw_1+x^{k+1}yw_1+x^ky(y*w_1),
\endalign$$
whose difference, assuming  the induction hypothesis, is
$$
xyx^{k-1}yw_1+\dots+x^ky^2w_1-x^{k+1}yw_1+x^ky(\bar D- D)(w_1).
\tag{11}$$
On the other hand, since $\bar D- D$ is a derivation, we have
$$
(\bar D- D)(x^kyw_1)=xyx^{k-1}yw_1+\dots+x^ky^2w_1-x^{k+1}yw_1+
x^ky(\bar D- D)(w_1) ,
$$
which agrees with (11).
\qed\enddemo
\par
\subheading{5. Action of the Hopf algebra $\QS$ on $\HH$}  In this
section we put a Hopf algebra structure on $\HH^1\cong\QS$, and define
an action of this Hopf algebra on $\QA$ that is related to the
results of the preceding sections.  In particular, we state a previous
result of the second author in terms of this action, and give a proof
of the result of Ihara and Kaneko stated in \S2.  All the algebraic 
definitions needed can be found in \cite{\K}.
\par
Clearly the
algebra $\QS$ of quasi-symmetric functions contains the algebra
Sym of symmetric functions.  In fact, the isomorphism $\phi:\HH^1\to\QS$
of the preceding section takes $z_n$ to the power-sum symmetric function 
$p_n=\sum_it_i^n$, and
$z_1^n$ to the elementary symmetric function $e_n=\sum_{i_1<\dots<i_n}t_{i_1}
\cdots t_{i_n}$.
Further, $\QS$ can be given a Hopf algebra structure that extends the usual
Hopf algebra structure on Sym (see \cite{\Gei,\Eh}); the primitives are
the power-sum symmetric functions $p_n$.
If $\De:\HH^1\to\HH^1\otimes\HH^1$ is the adjoint of the concatenation
product on the generators $z_i$, i.e.
$$
\De(z_{i_1}z_{i_2}\cdots z_{i_k})=
\sum_{j=0}^k z_{i_1}\cdots z_{i_j}\otimes z_{i_{j+1}}\cdots z_{i_k},
$$
then $(\HH^1,*,\De)$ is a Hopf algebra and
the map $\phi:\HH^1\to\QS$ is an isomorphism of Hopf algebras.
Henceforth we shall identify $(\HH^1,*,\De)$ with $\QS$ via
$\phi$; so $z_n$ is the $n$th power-sum symmetric function and so forth.
\par
Now define a $\bold Q$-linear map $\cdot:\HH^1\otimes\HH\to\HH$ as follows.  
Let $1\cdot w=w$ for any word $w$ of $\HH$.  For a nonempty word $u$ of 
$\HH^1$, let $u\cdot x = 0,$
$$
u\cdot y = \cases x^ky,&u=z_k,\\0,&\text{otherwise,}\endcases
$$
and
$$
u\cdot w_1w_2=\sum_u (u'\cdot w_1)(u''\cdot w_2)
\tag{12}$$
for words $w_1,w_2$ of $\HH$, where $\De(u)=\sum_uu'\otimes u''$.
Then the coassociativity of $\De$ 
implies that $u\cdot w\in\HH$ is well-defined for any words $u$ of 
$\HH^1$ and $w$ of $\HH$.
The following fact is immediate from the definitions.
\proclaim{Proposition 5.1} For words $w$ of $\HH$, $z_1\cdot w=D(w)$;
more generally, the linear map $D_n:\HH\to\HH$ defined by $D_n(w)=z_n\cdot w$
is the derivation sending $x$ to 0 and $y$ to $x^ny$.
\endproclaim
The map $\cdot$ is related to the multiplication $*$ of the
previous section as follows.
\proclaim{Lemma 5.2} For words $u\in\HH^1$ and $w\in\HH$,
$u\cdot w$ is the sum of terms in $u*w$ of length $\ell(w)$.
\endproclaim
\demo{Proof} We proceed by induction on $\ell(w)$.  If $\ell(w)\le 1$
then $w$ is either a power of $x$ or of the form $x^pyx^q$, and the
conclusion is clear from the definition.  Now suppose the conclusion
is true if $\ell(w)<n$ and let $w$ be a word of length $n$.  Writing
$u=z_i u_1$ and $w=x^{p-1}yw_1$, we have (from axiom (H3) above)
$$
u*w=z_i(u_1*w)+x^{p-1}y(u*w_1)+x^{i+p-1}y(u_1*w_1) .
$$
Note that only the last two terms can contribute words of length
$\ell(w)$.  Since $\ell(w_1)<n$, we have by the induction hypothesis
$$
\text{sum of terms of length $\ell(w)$ in }u*w=x^{p-1}y(u\cdot w_1)
+x^{i+p-1}y(u_1\cdot w_1) .
$$
But applying equation (12) to $u\cdot w=z_iu_1\cdot x^{p-1}yw_1$ gives
$$\align
u\cdot w&=(1\cdot x^{p-1}y)(u\cdot w_1)+(z_i\cdot x^{p-1}y)(u_1\cdot w_1)\\
&=x^{p-1}y(u\cdot w_1)+x^{i+p-1}y(u_1\cdot w_1) .
\qed\endalign$$
\enddemo
\proclaim{Theorem 5.3} The map $\cdot:\HH^1\otimes\HH\to\HH$ is an action
of the algebra $\QS\cong\HH^1$ on $\HH$, and in 
fact makes $\QA$ a $\QS$-module algebra.
\endproclaim
\demo{Proof} It suffices to show that $u\cdot(v\cdot w)=(u*v)\cdot w$
for words $u,v$ of $\HH^1$ and $w$ of $\HH$.  But by the lemma, both 
sides are just the sum of words of length $\ell(w)$ in $u*(v*w)=(u*v)*w$.
\qed\enddemo
\par
There is a relation between the action and the cyclic derivation $C$ of \S2.
\proclaim{Proposition 5.4} For positive integers $n,m$, $C(x^ny^m)=z_n\cdot
xy^m$.
\endproclaim
\demo{Proof} Using identity (5), we have
$$
(C(y^m),f)=xy^{m-1}fy+xy^{m-2}fy^2+\dots+xfy^m ,
$$
for any $f\in\HH$, so
$$
C(x^ny^m)=(C(x^ny^m),1)=(C(y^m),x^n)=xy^{m-1}x^ny+xy^{m-2}x^ny^2+\dots+
x^{n+1}y^m ,
$$
and the conclusion follows.
\qed\enddemo
This has the following corollary.
\proclaim{Corollary 5.5} For all $n,m\ge 1$, $\zt(z_n\cdot xy^m)=
\zt(z_m\cdot xy^n)$.
\endproclaim
\demo{Proof} By the preceding result, $z_n\cdot xy^m=C(x^ny^m)$; on the
other hand, $\tau\bar C(x^ny^m)=C(\tau(x^ny^m))=C(x^my^n)=z_m\cdot xy^n$,
and the conclusion follows from Theorems 2.3 and 4.1.
\qed\enddemo
A result of the second author \cite{\Oh} can be formulated in terms
of the $\QS$-action on $\QA$ as follows.  Let $h_n$ denote the complete
symmetric function of degree $n$, i.e. the sum of all monomials in
the $z_i$ of weight $n$.
\proclaim{Theorem 5.6} For all integers $n\ge 0$ and words
$w$ of $\HH^0$, $\zt(h_n\cdot\tau (w))=\zt(h_n\cdot w)$.
\endproclaim
\par
In view of Proposition 5.1 (and the $\tau$-invariance of $\zt$), 
Theorem 2.1 is the case $n=1$ of this theorem.  In fact, Ihara and
Kaneko proved the theorem stated in \S2 by showing it equivalent to
Theorem 5.6.  The argument that follows is based on their proof, but 
has been recast in terms of the $\QS$-action.
\par
If we let $\HH[[t]]$ be the ring of
formal power series in $t$ with coefficients in $\HH$, then the action
of $\QS$ on $\HH$ extends to an action of $\QS[[t]]$ on
$\HH[[t]]$.  Since $p_n=z_n$ acts on $\HH$ as the derivation $D_n$ of
Proposition 5.1, the operator
$$
\si_t=\exp\left(\sum_{n=1}^\infty \frac{D_n t^n}{n}\right)
$$
is an automorphism of $\HH[[t]]$ by the following result.
\proclaim{Lemma 5.7} Suppose $\de=t\de_1+t^2\de_2+\cdots$, where
each $\de_i$ is derivation of $\HH$.  Then
$$
\exp(\de)=\id+t\de_1+t^2\left(\frac{\de_1^2}{2}+\de_2\right)+\cdots
$$ 
is an automorphism of $\HH[[t]]$.
\endproclaim
\demo{Proof} First, note that $\de$ is a derivation of $\HH[[t]]$:  given
$u=u_0+tu_1+t^2u_2+\cdots$ and $v=v_0+tv_1+t^2v_2+\cdots$ in $\HH[[t]]$,
the coefficient of $t^n$ in $\de(uv)$ is 
$$
\sum_{p+q+r=n}\de_p(u_qv_r)=\sum_{p+q+r=n}(\de_p(u_q)v_r+u_q\de_p(v_r));
$$
but this is also the coefficient of $t^n$ in  $\de(u)v+u\de(v)$.  It then
follows that $\exp(\de)$ is an automorphism of $\HH[[t]]$, since for
$u,v\in\HH[[t]]$ we have
$$\multline
\exp(\de)(uv)=\sum_{n\ge 0}\frac{\de^n}{n!}(uv)=
\sum_{n\ge 0}\frac1{n!}\sum_{i=0}^n\binom{n}{i}\de^i(u)\de^{n-i}(v)=\\
\sum_{n\ge 0}\sum_{i+j=n}\frac{\de^i(u)}{i!}\frac{\de^j(v)}{j!}=
\exp(\de)(u)\exp(\de)(v). \qed
\endmultline$$
\enddemo
Now $H(t)=1+h_1t+h_2t^2+\cdots\in\QS[[t]]$, and from the well-known identity
$$
\frac{d}{dt}\log H(t)=\frac{H'(t)}{H(t)}=\sum_{n=1}^\infty p_n t^{n-1}
$$
(see, e.g., \cite{\M}) it follows that $\si_t(u)=H(t)\cdot u$ for 
$u\in\HH[[t]]$.  
Setting $\bar\si_t=\tau\si_t\tau$, we can restate Theorem 5.6 as saying
that $\zt(\bar\si_t(w)-\si_t(w))=0$ for any word $w$ of $\HH^0$, or
equivalently (since $\si_t(\HH^0)\subset\HH^0[[t]]$)
$$
\bar\si_t\si_t^{-1}(u)-u\in\ker\zt\quad\text{for all $u\in\HH^0[[t]]$}.
\tag{13}$$
Also, since
$H(t)^{-1}=E(-t)$, where $E(t)=1+e_1t+e_2t^2+\cdots=1+yt+y^2t^2+\cdots$,
we have $\si_t^{-1}(u)=E(-t)\cdot u$.  
\proclaim{Lemma 5.8} $\Phi=\bar\si_t\si_t^{-1}$ is uniquely 
characterized (among automorphisms of $\HH[[t]]$ that fix $t$) by the 
properties
\roster\item"{i.}"
$\Phi(x)=x(1-ty)^{-1}$
\item"{ii.}"
$\Phi(x+y)=x+y.$
\endroster
\endproclaim
\demo{Proof} To characterize an automorphism $\Phi$ of $\HH[[t]]$ that fixes
$t$, it is evidently enough to know where $\Phi$ sends $x$ and $y$; property
(i) gives $\Phi(x)$, and then property (ii) gives $\Phi(y)=x+y-\Phi(x)$.
To see that $\Phi=\bar\si_t\si_t^{-1}$ satisfies these properties, first note
that
$$\multline
\bar\si_t\si_t^{-1}(x)=\bar\si_t(E(-t)\cdot x)=
\bar\si_t(x)=\tau\si_t(y)=\tau(H(t)\cdot y)=\\
\tau(y+th_1\cdot y+t^2h_2\cdot y+\cdots)=
\tau(y+txy+t^2x^2y+\cdots)=x+txy+t^2xy^2+\cdots ,
\endmultline$$
and then do a similar calculation to show that $\bar\si_t\si_t^{-1}(y)=
y-txy(1-ty)^{-1}$.
\qed\enddemo
\par
Now consider the derivation 
$$
\partial_t=\sum_{n=1}^\infty t^n\frac{\partial_n}{n}
$$
of $\HH[[t]]$; by Lemma 5.7, $\exp(\partial_t)$ is an 
automorphism of $\HH[[t]]$.  To show that
$\partial_n(w)\in\ker\zt$ for all $n\ge 1$ and $w\in\HH^0$ is equivalent 
to (13) (and thus to Theorem 5.6), it suffices to prove the following.
\proclaim{Theorem 5.9} $\exp(\partial_t)=\bar\si_t\si_t^{-1}$.
\endproclaim
\demo{Proof} We use Lemma 5.8.  Since the derivations $\partial_n$
all take $z=x+y$ to 0, it is evident that $\exp(\partial_t)$ satisfies
property (ii).  To show $\exp(\partial_t)(x)=x(1-ty)^{-1}$, set
$$
G(s)=\exp(s\partial_t)(x)=\sum_{n=0}^{\infty}\partial_t^n(x)\frac{s^n}{n!}
\in\HH[[s,t]].
$$
Then $G(s)$ is the solution of the initial-value problem $G'(s)=
\partial_tG(s)$, $G(0)=x$.  We claim that
$$
G(s)=x\left(1-\frac{1-(1-tz)^s}{z}y\right)^{-1}
\tag{14}$$
since the right-hand side also satisfies these conditions; the conclusion 
then follows upon setting $s=1$.  
To verify the claim, let $U=(1-(1-tz)^s)/z$ and $V=\log(1-tz)/z$;
then the right-hand side of equation (14) is
$$
x(1-Uy)^{-1}=x(1+Uy+UyUy+\cdots)
$$
and the claim follows from the identities $U'(s)=UzV-V$, 
$\partial_tU=0$, $\partial_t(x)=-xVy$, and $\partial_t(y)=zVy-yVy$.
\qed\enddemo

\enddocument
\end